
\input amstex
\document
\magnification=\magstep1
\baselineskip 19.6pt
\overfullrule=0.pt

\par \hskip 2pc {\bf OPEN PROBLEMS ON OPEN ALGEBRAIC VARIETIES}

\par \vskip 1pc

\par \hskip 8pc (MONTREAL 1994 PROBLEMS)

\par \vskip 3pc

 A lively session held at the end of a conference on
Open Algebraic Varieties organized at the Centre de Recherches en
Math\'ematiques in December 1994 produced a list of open problems
that the participants would like to make available to the
mathematical community. Thanks are due to the contributors, to
D.-Q. Zhang, who undertook the initial collecting of the problems,
and to M. Zaidenberg, who was the guiding spirit behind the
efforts to prepare the collection for electronic distribution.

\par \vskip 3pc

 \par \hskip 15pc Peter Russell

\par \vskip 2pc

{\bf Contributors :}

\par \vskip 1pc

{\bf R.V. Gurjar} (gurjar$\@$tifrvax.tifr.res.in)

\par

{\bf Shulim Kaliman} (kaliman$\@$math.miami.edu)

\par

{\bf N. Mohan Kumar} (kumar$\@$artsci.wustl.edu)

\par

{\bf Masayoshi Miyanishi} (miyanisi$\@$math.sci.osaka-u.ac.jp)

\par

{\bf Peter Russell}  (russell$\@$Math.McGill.CA)

\par

{\bf Fumio Sakai} (fsakai$\@$rimath.saitama-u.ac.jp)

\par

{\bf David Wright} (wright$\@$einstein.wustl.edu)

\par

{\bf Mikhail Zaidenberg} (zaidenbe$\@$fourier.grenet.fr)

\par \vskip 5pc

{\bf R.V. Gurjar}

\par \vskip 1pc

{\bf PROBLEM.}  Let  $f: {\bold C}^n \rightarrow V$
be a proper (that is, a finite) morphism from the affine n-space onto
a normal affine variety  $V.$
Show that  $V$  is contractible.


\par \vskip 1pc

I proved
long time ago that

(i) $V$ is simply-connected, and

(ii) All homology groups of
$V$ are finite and the top homology group $H_n(V;Z)=(0)$.

In particular, if $n=2$ then $V$ is contractible. In fact, Miyanishi and
independently
Shastri-Gurjar proved that in this case $V$ is isomorphic to a quotient of
${\bold C}^2$
by a finite group.
Later on Shravan Kumar made a very nice improvement on this. He proved that
if $n=3$ then also $V$ is contractible. More generally, he proved that the
group
$H_{n-1}(V;Z)$ is trivial. His proof uses Smith theory in a nice way. The proof
has appeared in [K].

\par \vskip 1pc

REFERENCE :

[K]  Shravan Kumar, {\it A generalization of the Conner conjecture and topology
of Stein spaces dominated by ${\bold C}^n$}, Topology  {\bf 25} (1986),
483--494

\par \vskip 3pc

{\bf Shulim Kaliman}

\par \vskip 1pc

PROBLEM ABOUT CLASSIFICATION OF POLYNOMIALS IN TWO
\par
VARIABLES WITH A  ${\bold C}^*$-FIBER

\par \vskip 1pc

Let $p(x,y)$ and $q(x,y)$ be polynomials
in two complex variables.
We shall say that these polynomials are equivalent
if there exist a polynomial automorphism
$\alpha$  of  ${\bold C}^2$  and an affine
automorphism  $\beta$  of  ${\bold C}$  for which
$p = \beta \circ q \circ \alpha$.

\par \vskip 1pc

{\bf PROBLEM.}  Find the list of non-equivalent
polynomials such that every primitive
polynomial  $p$  with a  ${\bold C}^*$-fiber is
equivalent to one of polynomials from this list.
\par \vskip 1pc


\par
The answer to this question is known if one of the
following additional conditions holds [K2]

\par \vskip 1pc

(i)  $p$  is a rational polynomial, i.e. the genus of
its generic fibers is 0; or

\par
(ii) there exists a contractible curve in  ${\bold C}^2$  which
does not meet the  ${\bold C}^*$-fiber of  $p$.

\par \vskip 1pc

It is worth mentioning that P. Russell constructed a polynomial
with a  ${\bold C}^*$-fiber which satisfies neither of these
conditions.

\par \vskip 1pc

 More generally, consider the set $S_R$ of polynomials which have a fiber
isomorphic to a given affine algebraic curve $R$.
It is natural to look for a
list $L_R$ of non-equivalent polynomials such that every polynomial from $S_R$
is
equivalent to one of the polynomials from the list $L_R$.  If such
a list exists we
shall say that there is a classification of polynomials with this
fiber $R$.
This problem is equivalent to the problem of classification of
all smooth polynomial embeddings of $R$ into
${\bold C}^{2}$ up to a polynomial automorphism.  The
Abhyankar-Moh-Suzuki theorem  [AM], [Su]
says that all smooth polynomial embeddings
of the complex line into ${\bold C}^{2}$ are equivalent to linear
embeddings.  Moreover, V. Lin and M. Zaidenberg
[LZ] obtained the classification
of polynomial injections of ${\bold C}$ into ${\bold C}^{2}$ (i.e.
they found a description of all polynomials whose zero fiber is
homeomorphic to ${\bold C}$).  Later W. Neumann and L. Rudolph [NR]
reproved these theorems and W. Neumann obtained the classification for all
polynomials whose zero fiber is diffeomorphic
to a once-punctured Riemann surface
of genus $\leq 2$ [N].

The papers [AM], [NR], and [N] use essentially
the following theorem [AS]:

\par \vskip 1pc

{\it If the zero fiber of a polynomial is  a once-punctured
Riemann surface, then every other fiber of this polynomial is
once punctured.}

\par \vskip 1pc

The Lin-Zaidenberg theorem is based on the following elegant fact.

\par \vskip 1pc

{\it If a polynomial has at most one degenerate fiber (and it is so
in the case of a contractible fiber) then the polynomial is isotrivial,
i.e. its generic fibers are pairwise isomorphic.}

\par \vskip 1pc

Isotrivial
polynomials form a narrow class and its
classification was obtained later
in [K1], [Z1], [Z2].

If $R$ has more than one puncture none
of the above approaches works. The number of punctures on the generic
fiber of the corresponding polynomial may be arbitrary and the
polynomial may have a second degenerate fiber.  This makes
the problem difficult; that is why we suggest above to consider
$R$ isomorphic to ${\bold C}^*$ which is the simplest case
of a twice punctured Riemann surface.










\par \vskip 1pc

 REFERENCES:

[AM] S.S. Abhyankar, T.T. Moh, {\sl Embeddings of the line in the
plane}, J. Reine Angew. Math. {\bf 276} (1975), 148-166

[AS] S.S. Abhyankar, B. Singh, {\sl
Embeddings of certain curves in the affine plane}, Amer. J.
Math. {\bf 100} (1978), 99-175

[K1] Sh. Kaliman, {\sl Polynomials on $\bold C^2$ with isomorphic
generic fibers}, Soviet Math. Dokl. {\bf 33} (1986), 600--603

[K2] Sh. Kaliman, {\sl  Rational polynomial with a
$\bold C^*$-fiber},
Pacific J. Math. (to appear)

[LZ] V. Lin, M. Zaidenberg, {\sl An irreducible simply connected curve
in ${\bold C}^{2}$ is equivalent to a quasihomogeneous curve}, (English
translation) Soviet Math. Dokl. {\bf 28} (1983), 200-204

[NR] W.D. Neumann, K. Rudolph, {\sl  Unfoldings in knot theory},
Ann. Math. {\bf 278} (1987), 409-439 and Corrigendum:
{\sl Unfoldings in knot theory}, ibid {\bf 282},
349- 351

[N] W.D. Neumann, {\sl Complex algebraic plane curves via their link at
infinity}, Invent. Math. {\bf 98} (1989), 445-489

[Su] M. Suzuki, {\sl Propri\'{e}t\'es topologiques des
polyn\^{o}mes de
deux variables complexes et automorphismes alg\'{e}brigues de l'espace
${\bold C}^2$ }, J. Math. Soc. Japan {\bf 26} (1974), 241-257

[Z1] M. Zaidenberg, {\sl Ramanujam surfaces and exotic algebraic
structures on ${\bold C}^n$}, Dokl. AN SSSR 314(1990), 1303--1307, English
translation in Soviet Math. Doklady {\bf 42} (1991), 636--640

[Z2] M. Zaidenberg, {\sl On Ramanujam surfaces, ${\bold C}^{**}$-families
and exotic algebraic structures on ${\bold C}^n$, $n \ge 3$},
Trudy Moscow Math. Soc. {\bf 55} (1994), 3--72 (Russian; English transl.
to appear)

\par \vskip 3pc

{\bf N. Mohan Kumar}

\par \vskip 1pc

Let  $A$  be an Artin local ring over an algebraically closed
field $k$ of characteristic 0, with $A$ finite dimensional
over $k$ (as a vector space). For any finitely generated
module  $M$  over  $A$, let  $D(M)$  denote the module of
$k$-derivations from  $A$  to  $M$.  Also let  $l$  stand
for length of a module.

\par \vskip 1pc

{\bf QUESTION 1.}  Is it true that  $l(\Omega_A^1) \geq l(A)-1$  and
equality if and only if  $A$  is isomorphic to  $k[x]/x^n$ ?

\par \vskip 1pc

{\bf QUESTION 2.}  More generally, is it true that
$l(D(M)) \geq l(M)-l(s(M))$  where  $s(M)$  denotes
the socle of  $M$  and equality if and only if  $A$
is isomorphic to  $k[x]/x^n$  and  $M$  free ?

\par \vskip 3pc

{\bf Masayoshi Miyanishi}

\par \vskip 1pc

VECTOR FIELDS ON FACTORIAL SCHEMES

\def\Z{\bold Z}
\def\A{\bold A}
\def\BP{\bold P}

\def\Spec{{\text{\rm Spec}}\:}
\def\Proj{{\text{\rm Proj}}\:}
\def\lto{\longrightarrow}
\def\dlto{\:\cdots\!\!\to}
\def\trdeg{{\text{\rm tr.deg}}\:}

\def\dps{\displaystyle}

\par \vskip 1pc
Let $K$ be a field of characteristic zero, let $R$ be
a noetherian $K$-algebra domain and let $L$ be the
quotient field of $R$. Let $\delta$ be a $K$-derivation
on $R$. We call an ideal $I$ of $R$ a
{\bf $\delta$-integral ideal} if $\delta(I) \subseteq I$.
If a principal ideal $I = fR$ is a  $\delta$-integral ideal,
$f$ is called a {\it $\delta$-integral element}.
When we restrict ourselves to the case $R = K[x,y]$,
we call $f$ a {\it $\delta$-integral curve}.
If $f$ is a $\delta$-integral element, we can
write $\delta(f) = f\chi(f)$ with $\chi(f) \in R$.
An element $t$ of $R$  is called a
{\it $\delta$-integral factor} if there exists a
$\delta$-integral element $f$ such that $t = \chi(f)$.
The set of all $\delta$-integral factors in $R$ is
denoted by $X_{\delta}(R)$  or simply by $X_{\delta}$.
We note that an invertible element of $R$ is a
$\delta$-integral element.

\par \vskip 1pc

LEMMA 1.\ {\it {\rm (1)}\ \ Let $t$ be a $\delta$-integral
factor, and let $A_t$  be the set of $\delta$-integral
elements in $R$ with the same  $\delta$-integral factor $t$.
Then $A_t$ is a $K$-vector space.

\par

{\rm (2)}\ \ $X_{\delta}$ is an abelian monoid under
the addition of $R$, i.e., $X_{\delta}$ is closed
under the addition and contains the zero element.

\par

{\rm (3)}\ \ Let $A$ be the subalgebra of $R$ generated
by all  $\delta$-integral elements over $K$.
Then $\dps{A = \sum_{t \in X_{\delta}}A_t}$ and
$A_s\cdot A_t \subseteq A_{s+t}$, while
$\dps{\sum_{t\in X_{\delta}}A_t}$ is not necessarily
a direct sum. We call $A$ the {$\delta$-integral ring} of $R$.}

\par \vskip 1pc

Let $L = Q(R)$ be the quotient field of $R$. Then the
$K$-derivation  $\delta$ on $R$ is naturally extended
to a $K$-derivation on $L$.  An element $\xi$ of $L$ is,
by definition, a $\delta$-integral element  in $L$ if
$\chi(\xi) = \delta(\xi)/\xi \in R$, and $\chi(\xi)$ is
then a  $\delta$-integral factor. The set of all
$\delta$-integral factors in $L$  is an abelian group
under the addition, which we denote by $\widetilde{X_{\delta}}$.
Clearly, $X_\delta \subseteq \widetilde{X_{\delta}}$ as
monoids. The next lemma shows that
$\widetilde{X_{\delta}} = X_{\delta}-X_{\delta}$,
i.e., $\widetilde{X_{\delta}}$ is an abelian group
generated by $X_{\delta}$ provided $R$ is a factorial
domain (= a unique factorization domain).

\par \vskip 1pc

LEMMA 2.  \ {\it Assume that $R$ is a factorial domain.
Then we have:

\par

{\rm (1)}\ \ Let $\xi$ be a $\delta$-integral element
in $L$ and write  $\xi = fg^{-1}$ with mutually prime
elements $f, g$ of $R$. Then $f$ and $g$ are $\delta$-integral
elements in $R$, and $\chi(\xi) = \chi(f) -\chi(g)$.

\par

{\rm (2)}\ \ Let $f$ be a $\delta$-integral element in $R$.
Then any prime factor as well as any divisor of $f$ is
a $\delta$-integral element in $R$.
Furthermore, the $\delta$-integral ring $A$ is generated
over $K$ by invertible elements of $R$ and $\delta$-integral
elements which are prime elements of $R$.}

\par \vskip 1pc

We then pose the following:

\par \vskip 1pc

{\bf QUESTION 1.} \ {\it With the notations and assumptions
as above, is $\widetilde{X_{\delta}}$  a finitely
generated abelian group provided $R$ is finitely
generated over $K$ ?}

\par \vskip 1pc

Later, we shall present one result which asserts
that $\widetilde{X}_{\delta}$  is finitely generated.

\par \vskip 1pc

DEFINITION. \ Let $X$ be an abelian monoid.

\par

{\rm (1)}\ \ $X$ is {\it positive} if $X$ contains
no abelian subgroups other than $(0)$.

\par

{\rm (2)}\ \ $X$ is {\it finitely generated} if
$\widetilde{X}:= X-X$ is a finitely generated abelian group.

\par

{\rm (3)}\ \ $X$ is {\it good} if $X$ is positive
and if there exist elements $t_1, \cdots, t_r$ of $X$
such that $X = \Z_+t_1+\cdots +\Z_+t_r$ and
$\widetilde{X}$ is a free abelian group with
free basis  $t_1, \cdots, t_r$, where $\Z_+$ is
the set of non-negative integers. We then write
$X = \Z_+t_1\oplus \cdots \oplus \Z_+t_r$
and call $r$ the rank of $X$.

\par \vskip 1pc

We shall consider the subring
$A_0 = \{x \in R; \delta(x) = 0\}$ of $R$
and the subfield $L_0 = \{\xi \in L; \delta(\xi) = 0\}$
of $L$. We say that $A_0$ is an {\it inert} subring
of $R$ if $a \in A_0$ and $a = bc$ with
$b, c \in R$ implies $b, c \in A_0$.

\par \vskip 1pc

LEMMA 3. \ {\it Let $R$ be a factorial domain. Then we have:

\par

{\rm (1)}\ \ $X_{\delta}$ is positive if and only
if $A_0$ is an inert subring of $R$.

\par

{\rm (2)}\ \ $L_0$ is algebraically closed in $L$.
If $A_0$ is an inert subring of $R$ then $Q(A_0)$
is algebraically closed in $L$.

\par

{\rm (3)}\ \ Every element $\xi$ of $L_0$ is written
as $\xi = b/a$, where $a$ and $b$ are $\delta$-integral
elements in $R$ with the same  $\delta$-integral factor.
Conversely, if $a, b \in A_t$ with $t \in X_{\delta}$
and $a \neq 0$ then $b/a \in L_0$.}

\par \vskip 1pc

We have the following result.

\par \vskip 1pc

THEOREM 4. \ {\it Let $R$ be a noetherian $K$-algebra
domain and let $\delta$ be a $K$-derivation on $R$.
Assume that $R$ is a factorial domain and the monoid
$X_\delta$ of $\delta$-integral factors is good.
Write $X_{\delta} = \Z_+t_1\oplus \cdots \oplus \Z_+t_r$
with $t_1, \cdots, t_r \in R$. Let $f_i\ \ (1 \leq i \leq r)$
be a $\delta$-integral element such that $\chi(f_i) = t_i$.
Then the following assertions hold.

\par

{\rm (1)}\ \ We may assume that $f_1, \cdots, f_r$
are prime elements of $R$.

\par

{\rm (2)}\ \ $R^* \subseteq A_0$ and $A_0$ is an inert subring of $R$.

\par

{\rm (3)}\ \ For the $\delta$-integral ring $A$ of $R$,
 we have $A\otimes_{A_0}L_0 = L_0[f_1, \cdots, f_r]$.
If $f_1, \cdots, f_r$ are algebraically independent over
$L_0$, $A = \dps{\sum_{t\in X_{\delta}}A_t}$ is a graded ring.
Namely, the decomposition $\dps{\sum_{t\in X_{\delta}}A_t}$
is a direct sum.

\par

{\rm (4)}\ \ $L_0 = Q(A_0)$ if and only if
$A = A_0[f_1, \cdots, f_r]$.
Furthermore, if $L_0 = Q(A_0)$ and
$\dps{A = \sum_{t\in X_{\delta}}A_t}$ is a graded ring,
then $f_1, \cdots, f_r$ are algebraically independent over $L_0$.}

\par \vskip 1pc

LEMMA 5. \ {\it Assume that $L_0 = Q(A_0), X_\delta$
is positive and $\widetilde{X_{\delta}}$ is finitely generated.
Then $X_{\delta}$ is good.}

\par \vskip 1pc

DEFINITION. \ We say that $\delta$ is {\it locally nilpotent}
if, for each $x \in R, \delta^n(x) = 0 \ (n \gg 0)$.
Let $T$ be an indeterminate. Define a mapping
$\varphi : R \lto R[T]$ by
$\varphi (x) = \dps{\sum_{i\geq 0}\frac{1}{i !}}\delta^i(x)T^i$.
Then $\varphi$ is a homomorphism of $K$-algebras.

\par \vskip 1pc

LEMMA 6. \ {\it Suppose $\delta$ is locally nilpotent.
Then $X_{\delta} = (0), A = A_0$  and $L_0 = Q(A_0)$.}

\par \vskip 1.5pc

Now, we assume that $R$ is a factorial domain of
dimension $2$ which is finitely generated over $K$.
For a $K$-derivation $\delta$ of $R$, the set
$\delta(R) = \{\delta(x); x \in R\}$ generates an
ideal $R\delta(R)$.  The {\it divisorial part}
$(\delta)$ of $\delta$ is the greatest principal
ideal of $R$ which divides $R\delta(R)$. Set
$(\delta) = dR$ with $d \in R$.
Then $\delta' = d^{-1}\delta$ is a $K$-derivation of
$R$ such that the ideal  $R\delta'(R)$ has height $2$;
we then say that $\delta'$ has  {\it no divisorial part}.
Let $V =$ Spec $(R)$. The subset  $V(R\delta(R))_{\text{\rm red}}$
is called the {\it zero set} of $\delta$, which
we denotes by $Z(\delta)$. If $\delta$ has no
divisorial part, $Z(\delta)$ is a finite set.

\par \vskip 1pc

LEMMA 7. \ {\it Let $K \subset K_1 \subset L:= Q(R)$
be a subfield such that  $\trdeg_KK_1 = 1$.
Then there exists a non-trivial $K$-derivation $\delta$ of
$R$ determined uniquely up to an invertible element
of $R$ such that $\delta$ has no divisorial part and
$L_0 := \{\xi \in L;\delta(\xi) = 0\} \supset K_1$.}

\par \vskip 1pc

Let $\varphi : V = \Spec R \dlto C$ be a rational
mapping onto a smooth algebraic curve $C$.
Then $\varphi$ may not be defined in a finite set
$\Sigma$ of $V$. Namely,
$\varphi^0 := \varphi|_{V-\Sigma} : V-\Sigma \lto C$
is a morphism. For $P \in C$, the schematic
closure of the fiber $(\varphi^0)^{-1}(P)$ is called
the fiber of $\varphi$ over $P$ and denoted by
$\varphi^{-1}(P)$. We say that a rational mapping
$\varphi : V \dlto C$ is an {\it irreducible pencil}
parametrized by $C$ if $\varphi^0$ is surjective
and general fibers of $\varphi$ are irreducible and
reduced. A rational mapping $\varphi : V \dlto C$ is
equivalently defined by giving a subfield $K(C)$ of
$L := Q(R)$. Then $\varphi$ is an irreducible pencil
if and only if $K(C)$ is algebraically closed in $L$.
Given a rational mapping $\varphi : V \dlto C$, the
Stein factorization
$\varphi : V \overset{\widetilde{\varphi}}\to{\dlto}$
$\widetilde{C} \overset{\sigma}\to{\lto} C$
gives an irreducible pencil
$\widetilde{\varphi}: V \dlto \widetilde{C}$.
An irreducible pencil $\varphi : V \dlto C$ is called a
{\it fibration} if $\varphi$ is a morphism.

\par

Let $\delta$ be a $K$-derivation of $R$. We say that
$\delta$ is {\it composed} of a pencil
$\varphi : V \dlto C$ if $\delta(I_F) \subseteq I_F$
for every general fiber $F$ of $\varphi$,
where $I_F$ signifies the defining ideal of $F$.
We also say that $\delta$ is {\it of fibered type}
if $\delta$ is composed of some pencil.

\par \vskip 1pc

LEMMA 8 \ {\it Let the notations be the same as above.
Assume that $R^* = K^*$, where $R^*$ is the group of
invertible elements of $R$ and that $\delta$ is composed
of an irreducible pencil $\varphi : V \dlto C$.
Then the following assertions hold:

\par

{\rm (1)}\ \ $C$ is either $\A^1$ or $\BP^1$.

\par

{\rm (2)}\ \ In case $C = \A^1$, there exists an
invertible element $f$ of $R$ with $\delta(f) = 0$
such that $C = \Spec K[f]$ and $\varphi$ is a morphism
associated with the inclusion $K[f] \hookrightarrow R$.
We then have $A_0 \supset K$.

\par

{\rm (3)}\ \ In case $C = \BP^1$, there exists two
irreducible elements $f, g$ of $R$ such that
$\chi(f) = \chi(g), \gcd(f,g) = 1$ and
$\varphi : V \dlto \BP^1$ is the natural extension of
a fibration $D(g) \lto \Spec K[f/g] = \A^1$, where
$\BP^1 = \Proj K[f,g]$.
We then have $L_0 \supset A_0 = K$.}

\par \vskip 1pc

LEMMA 9. \ {\it Assume that $R^* = K^*$. Let $f$ be
an irreducible element of $R$ and let $\varphi : V \lto \A^1$
be a fibration defined by the inclusion $K[f]
\hookrightarrow R$. Let $\delta$ be a $K$-derivation of $R$
such that $L_0 \supseteq K(f)$ and $\delta$ has no
divisorial parts (cf. Lemma 7).
Then the following assertions hold:

\par

{\rm (1)}\ \ $\delta$ is composed of the fibration $\varphi$.

\par

{\rm (2)}\ \ The following conditions are equivalent to each other:

\par

{\rm (i)}  $\trdeg_KQ(A) = 1$.

\par

{\rm (ii)} $A = A_0 = K[f].$

\par

{\rm (iii)}  $X_\delta = (0)$.

\par

{\rm (iv)}  $X_\delta$ is positive.

\par

{\rm (v)}  Every fiber of $\varphi$ is irreducible.

\par

{\rm (3)}\ \ Assume that $R$ is a polynomial ring  $K[x,y]$. Then
$\delta \,\, = \dps{d^{-1}\left( \frac{\partial f}{\partial y}
\frac{\partial}{\partial x}-\frac{\partial f}{\partial x}
\frac{\partial}{\partial y}\right)},$
where  $d = \gcd \dps{\left(\frac{\partial f}{\partial x},
\frac{\partial f}{\partial y}\right)}$.
Furthermore, if $d = 1$, the following conditions
are equivalent to each other:

\par

{\rm (i)} $\delta$ has no zeroes, i.e., $Z(\delta) = \emptyset$.

\par

{\rm (ii)}  Every fiber of $\varphi$ is smooth.

\par \vskip 1pc

If $X_\delta = (0)$ then $d = 1$.}

\par \vskip 1pc

LEMMA 10. \ {\it Assume that $R^* = K^*$. Let $f, g$
be two (distinct) irreducible elements of $R$ such that
$1 \not\in Kf+Kg$ and $fR+gR = R$, and let
$\varphi : V \dlto \BP^1$ be a rational mapping defined
by $P \longmapsto (f(P):g(P))$.  Let $\delta$ be a non-trivial
$K$-derivation of $R$ such that $\delta$ has
no divisorial parts and $L_0 \supset K(f/g)$ (cf. Lemma 7).
Then the following assertions hold:

\par

{\rm (1)}\ \ $\delta$ is composed of the pencil $\varphi$.

{\rm (2)}\ \ $L_0 = K(f/g)$ and $A_0 = K$.

{\rm (3)}\ \ In the case where $R$ is a polynomial ring $K[x,y]$,
$\delta$  is determined, up to a constant multiple, as
$$
\delta = d^{-1}\left(g\left(\frac{\partial f}{\partial y}\frac{\partial}
{\partial x}-\frac{\partial f}{\partial x}\frac{\partial}{\partial y}\right)
-f\left(\frac{\partial g}{\partial y}\frac{\partial}{\partial x}-
\frac{\partial g}{\partial x}\frac{\partial}{\partial y}\right)\right),
$$
where
$$
d = \gcd\left(g\frac{\partial f}{\partial y}-f\frac{\partial g}{\partial y},
g\frac{\partial f}{\partial x}-f\frac{\partial g}{\partial x}\right).
$$
Furthermore, if $d = 1$ and $J(f,g) \in K^*$
then $Z(\delta) = \emptyset$, where $J(f,g)$ is the Jacobian
of $f, g$ with respect to $x, y$.}

\par \vskip 1pc

We say that a $K$-derivation $\delta$ of $R$ is locally
nilpotent along an irreducible pencil fibration
$\varphi : V \dlto C$ if $\delta$ is composed of $\varphi$
and if the restriction $\delta_F$ of $\delta$ onto $F$ is a
locally nilpotent $K$-derivation of $R/I_F$, where $F$ is a
general fiber of $\varphi$ and $I_F$ is the defining ideal of $F$.

\par \vskip 1pc

THEOREM 11. \ {\it Let $R$ be a factorial domain of
dimension two which is finitely generated over $K$ and
let $\delta$ be a non-trivial $K$-derivation. Assume that
$R^* = K^*$. Then the following conditions
are equivalent to each other:

\par

{\rm (1)}  $R = K[x,y]$, a polynomial ring, and $\delta(y) \in K[x]$.

\par

{\rm (2)}  $\delta$ is locally nilpotent.

\par

{\rm (3)}  There exists a fibration $\varphi : V \lto C$
such that $\delta$ is locally nilpotent along $\varphi$.}

\par \vskip 1pc

THEOREM 12. \ {\it Let $R = K[x,y]$ be a polynomial ring,
let $f$ be an irreducible polynomial in $R$ and let $\delta =
\dps{\frac{\partial f}{\partial y}\frac{\partial}{\partial x}-
\frac{\partial f}{\partial x}\frac{\partial}{\partial y}}$.
Assume that $\delta$ has no divisorial parts. Then $\delta$
is locally nilpotent if and only if $X_{\delta} = (0)$ and
there exists a curve $C$ which is defined by $g = 0$ with a
$\delta$-integral element $g$ of $R$ and isomorphic to $\A^1$.}

\par \vskip 1pc

The following result gives a partial answer to the question
on finite generation of $\widetilde{X}_\delta$.

\par \vskip 1pc

THEOREM 13. \ {\it Let $R$ be a factorial domain of
dimension two which is finitely gnerated over $K$ and
let $\delta$ be a non-trivial $K$-derivation of fibered
type on $R$. Assume that $R^* = K^*$.
Then $\widetilde{X}_\delta (= X_{\delta}-X_{\delta})$
is finitely generated.}

\par \vskip 1pc

If $A_0 \supset K$ or $L_0 \supset K$,
a derivation $\delta$ if of fibered type. We say that
a $K$-derivation $\delta$ is {\it of general type}
if $A_0 = L_0 =K$. Then our second question is the following:

\par \vskip 1pc

{\bf QUESTION 2.} \ {\it Classify all $K$-derivations of general
type on $R = K[x,y]$ and describe
them in terms of the $\delta$-integral ring $A$.}

\par \vskip 1.5pc

REFERENCE :

M. Miyanishi, {\sl Vector fields on fractional schemes},
to appear in J. Algebra.

\par \vskip 3pc

{\bf Peter Russell}

\par \vskip 1pc

{\bf QUESTION 1.}  Let  $A$  be a purely inseparable form of
$k^{[n]},$  the polynomial ring in  $n$  variables over the
field  $k,$  i.e. $A$  is isomorphic to  $K^{[n]}$  over
the perfect closure  $K$  of  $k.$
Is  $A$  trivial if  $A$  is $k$-rational, i.e. if the
quotient field  $L$  of  $A$  is isomorphic to  $k^{(n)},$ or,
more strongly, if  $A$  is birationally contained in  $k^{[n]}$ ?

\par \vskip 1pc

Not much is known about purely inseparable forms of $k^{[n]}$ except
under strong homogneity asumptions, and even for $n = 1$ our information
is quite incomplete [KMT] [KM] [R1]. Our question, easily answered
for $n = 1$, arose in efforts to extend the results of [BR] on the
behaviour of birational subrings of $k^{[2]}$ under basefield extension
to the purely inseparable case.

\par \vskip 1pc

{\bf QUESTION 2.}  Let  $(A, M)$  be a regular local ring and
$B = A[a/b]$  with  $a$  in  $M \setminus M^2$  and  $b$  in  $M^2$
and  $GCD(a, b) = 1.$
(These conditions are there to make  $B$  regular.)
Under what conditions is  $MB$  a complete intersection in  $B$ ?
(An obvious sufficient condition is that  $a$  can be
extended to a regular system of parameters  $a, u, ..., v$
with  $b$  in  $(u,...,v)A.$
If  ${\text dim}\,(A) = 2,$ this condition is also necessary.)

\par \vskip 1pc

Let $A = {\bold C}[x,y,z,t]$ with $x + x^2y + z^2 + t^3 = 0$. The problem came
out of efforts to decide whether the ideal $(x,z,t)$ is a complete intersection
in $A$. This is of interest since it is known that ${\text Spec}\,(A)$
is a smooth contractible threefold [KR] not isomorphic to ${\bold C}^{[3]}$
[ML].

\par \vskip 1pc

{\bf QUESTION 3.}  If the polynomial  $p$  in two variables over
${\bold C}$  has smooth, irreducible rational zero set,
is  $p$  a factor of a field generator
(or generically rational polynomial),
i.e. does there exist a polynomial  $q$  such that
the generic fibre of  $pq$  is rational ?

\par \vskip 1pc

To "effectively" describe all polynomials $p$ in two variables with
smooth rational zero set is certainly a formidable task. There is a
bit more hope for those $p$ that give a fibration of the plane with
rational generic fibre. (Such $p$ have been widely studied [MS] [S]
[R2] [R3]. See also the problem proposed by Sh. Kaliman and the
references given there.) It is known that every irreducible factor
$q$ of $p-c,\, c \in \bold C$, has smooth rational zero set, but it is not hard
to find examples where then the zeros of $q-c$ are non rational for
general $c \in \bold C$. We propose to investigate whether some converse
might be true.

\par \vskip 1pc

REFERENCES:

[BR]  S. Bhatwadekar and P. Russell, {\sl A note on geometric factoriality},
      1994 (to appear in Can. Math. Bull.)

[KM]  T. Kambayashi and M. Miyanishi, {\sl Forms of the affine line over a
      field}, Kinokuniya 1977, Tokyo

[KMT] T. Kambayashi, M. Miyanishi, M. Takeuchi, {\sl Unipotent algebraic
      groups}, Springer Lecture Notes in Mathematics {\bf 414} 1974

[KR]  M. Koras and P. Russell, {\sl Contractible threefolds and $\bold
C^*$-actions
      on ${\bold C}^3$} (in preparation)

[ML]  L. Makar-Limanov, {\sl On the hypersurface $x + x^2y + z^2 + t^3 = 0$},
      preprint, 1994

[MS]  M. Miyanishi and T. Sugie, {\sl Generically rational polynomials},
      Osaka J. Math. {\bf 17} (1980), 339-362

[R1]  P. Russell, {\sl Purely inseparable forms of the affine line and its
      additive group}, Pacific J. Math. {\bf 32} (1970), 527-539

[R2]  P. Russell, {\sl Field generators in two variables}, J. Math. Kyoto U.
      {\bf 15} (1975), 555-571

[R3]  P. Russell, {\sl Good and bad field generators}, J. Math. Kyoto U.
      {\bf 17} (1977), 319-331

[S]   H. Saito, {\sl Fonctions enti\`eres qui se reduisent a certains
      polynomes (II)}, Osaka J. Math. {\bf 14} (1977), 649-674

\par \vskip 3pc

{\bf Fumio Sakai}

\par \vskip 1pc

{\bf PROBLEM 1.}  Classify all rational and elliptic cuspidal
plane curves.

\par \vskip 1pc

{\bf PROBLEM 2.}  Find the maximal number of cusps among
all the rational cuspidal plane curves.

\par \vskip 3pc

{\bf David Wright}

\par \vskip 1pc


{\bf QUESTION.}  Can one prove there does not exist a
counterexample $(f,g)$ to the Jacobian Conjecture such that
$$f,g\in {\bold C}[Y,XY,X^2Y,X^3Y-X]$$
Note:  A positive answer to this question would be
a step toward proving the 2-dimensional Conjecture
in the case of smooth integral closure.

\par \vskip 3pc

{\bf Mikhail Zaidenberg}

\par \vskip 1pc



\S 1. ACYCLIC SURFACES

\par \vskip 1pc

{\bf 1.1.}  Let $X$ be a smooth contractible complex affine
algebraic surface.  If  $\overline{\kappa}(X) = 1$, where
$\overline{\kappa}$  denotes the logarithmic Kodaira dimension,
then there is only one simply connected curve $l_X$ in  $X$,
and this curve is isomorphic to  ${\bold C}$
( see [GuMi, GuPa, MiTs, Za 1(Addendum)]).
Let  $l_X = \{p = 0\}$, where  $p$  is an irreducible
regular function on  $X$.  Then all the other fibres
$F_c = p^{-1}(c),$ $\,c \neq 0$, of  $p$  are pairwise
isomorphic smooth once-punctured curves.  The image of
any morphism from a once-punctured curve  $\Gamma$
into  $X$  is known to be contained in a fibre of  $p$
(Sh. Kaliman, L. Makar--Limanov [KML 1]). Thus, $X$ containes,
up to an isomorphism, only two such curves, namely
$l_X = F_0$  and  $F_X := F_1$.

\par
Consider now a smooth contractible surface $X$ of
log-general type. It is known that there is no
simply-connected curve in $X$ (M. Zaidenberg [Za 1]; see also [GuMi, MiTs]).

\par \vskip 1pc

{\bf QUESTION.}  Does  $X$  contain any once-punctured curve ?
Does it contain only a finite number of such curves ?

\par \vskip 1pc

Of course, the same question has sense for
$\bold Q$-acyclic surfaces of log-general type, as well
as for embeddings of  ${\bold C}^*$  into  $X$
(the latter question has been proposed by M. Miyanishi; oral communication).

\par \vskip 1pc

{\bf 1.2.}  A smooth contractible algebraic surface $X$
is known to be affine (T. Fujita [Fu]).  Consider a proper
embedding  $X \hookrightarrow {\bold C}^N$  and the
plurisubharmonic function  $\varphi(x) = ||x||^2$  on  $X.$
Put  $X_R = \{\varphi < R^2\}.$  Then for all sufficiently
large $R$, \,\, $X_R$ is a strictly pseudoconvex domain
in  $X$  with a real-algebraic boundary. Furthermore,
$X_R$  is  diffeomorphic to  $X$  and hence contractible.
Fix such an  $R$.  Note that  $X_R$  is homeomorphic to
a bounded domain in  ${\bold R}^4.$

\par \vskip 1pc

{\bf QUESTION.}  Assume that  $X$  as above is not isomorphic to
${\bold C}^2$. Is it true that then  $X_R$  is not biholomorphic
to a bounded (strictly) pseudoconvex domain in  ${\bold C}^2 $ ?

\par \vskip 1pc

The boundary  $S_R = \partial X_R$  is a homology sphere
with a non-trivial perfect fundamental group (C. P. Ramanujam [Ra]).
So, more generally we may ask

\par \vskip 1pc

{\bf QUESTION.}  Let a homology sphere  $S$  be a boundary of
a strictly pseudoconvex domain in  ${\bold C}^2$.
Is it true that then  $S$  is diffeomorphic to  $S^3 $ ?

\par \vskip 1pc

{\bf 1.3.} {\bf RIGIDITY CONJECTURE} (H. Flenner-M. Zaidenberg [FlZa]).
Any smooth  $\bold Q$-acyclic surface  $X$  of log-general
type is rigid.

\par \vskip 1pc

This conjecture has been verified in a number of examples [FlZa].

\par \vskip 1pc

Consider a minimal SNC-completion  $(V,\,D)$  of  $X$, i.e.
$X = V - D$, where  $V$  is a smooth projective surface
and  $D = \sum_i D_i$  is a simple normal crossings divisor
in  $V$, which is minimal in this class.   The Euler characteristic of the
logarithmic tangent bundle of
$(V,\,D)$  can be expressed as follows:
$$\chi(\Theta_V\langle D \rangle) = K_V (K_V + D) =
10 - 3r - \sum_i D_i^2 \,\,,$$
where  $K_V$  is the canonical divisor of   $V$ [FlZa].
The Rigidity Conjecture induces the following one [FlZa]:

\par \vskip 1pc

{\bf 1.4.} {\bf CONJECTURE.}  For $X,\, V,\, D$  as above
$K_V (K_V + D) = 0.$

\par \vskip 1pc

Since  $D$  is a rational tree, this is equivalent to
the equality for the logarithmic Chern number
${\overline c}_1^2 (X) :=  (K_V + D)^2 = -2$. Together
with the evident equality for the topological Euler
characteristic  ${\overline c}_2 (X) = e(X) = 1$  this would
imply that  $3{\overline c}_2 (X) - {\overline c}_1^2
(X) = 5$, which was conjectured by T. tom Dieck [tD].

\par \vskip 1pc

{\bf 1.5.}  Let things be as above.  Denote by  $\Gamma_D$
the (unweighted) dual graph of  $D$ (which is a tree),
and by  $\breve \Gamma_D$  the (unweighted)
Eisenbud--Neumann diagram of  $D$, i.e. a tree which is
obtained from  $\Gamma_D$  by replacing every linear branch
between two neighbouring branching points of  $\Gamma_D$
by a single line, and every extremal linear branch
of $\Gamma_D$ by an arrowhead.

\par \vskip 1pc

{\bf 1.6.} {\bf QUESTION.}  Is the number of all Eisenbud--Neumann
diagrams $\breve \Gamma_D$ of the boundaries  $D$  of the minimal smooth
SNC--completions of the contractible
(resp. acyclic, resp. ${\bold Q}$-acyclic) surfaces of log--general type finite
?

\par \vskip 1pc

{\bf 1.7.} It is easily seen that for any non-constant polynomial $p \in {\bold
C}[x, y]$ there exists another one $q \in {\bold C}[x, y]$ such that the
regular mapping $F = (p,\,q)\,:\,{\bold C}^2 \to {\bold C}^2$ is proper. On the
other hand, if $X$ is a contractible surface of log--Kodaira dimension $1$ and
if $l_X = p^* (0),\, p \in {\bold C}[X]$, is the unique curve in $X$ isomorphic
to ${\bold C}$, then there is no $q \in {\bold C}[X]$ such that $F =
(p,\,q)\,:\,X \to {\bold C}^2$ would be a proper mapping [Za 1, Addendum].

\par \vskip 1pc

{\bf OUESTION.}  Let  $X$  be a smooth acyclic surface such that for any
non--constant regular function $p \in {\bold C}[X]$ there exists another one $q
\in {\bold C}[X]$ with the property that $F = (p,\,q)\,:\,X \to {\bold C}^2$ is
a proper morphism. Is it true that $X$ is isomorphic to  ${\bold C}^2$ ?

\par \vskip 1.5pc

\S 2.  EXOTIC ${\bold C}^n$'s

\par \vskip 1pc

{\bf 2.1.}  {\bf QUESTION.}  Let  $n > 3$.  Does there exist any
smooth contractible hypersurface of  ${\bold C}^n$  of
log-general type ?

\par \vskip 1pc

{\bf 2.2.}  {\bf QUESTION.}  Let $p \in {\bold C}[x_1 ,\dots, x_n]$
be such that the hypersurface  $X = \{p = 0 \}$  in
${\bold C}^n$  is contractible.  Is then  $p$  an isotrivial
polynomial, i.e. such that its generic fibres are pairwise
isomorphic hypersurfaces in ${\bold C}^n$ ?

\par See the problems of Sh. Kaliman above and the references therein on the
description of the isotrivial polynomials on ${\bold C}^2$. By the way, is it
possible to describe those polynomials on ${\bold C}^2$ which are {\it
projectively isotrivial} in the sense that their generic fibres all have
isomorphic smooth projective models? All generically rational polynomials (i.e.
{\it field generators}, cf. Question 3 of P. Russell) are in this class.

\par \vskip 1pc

{\bf 2.3.}  {\bf QUASIHOMOGENEITY CONJECTURE.}  If  $X$  as above
has an isolated singular point, then the polynomial $p$ is
equivalent, up to an automorphism of  ${\bold C}^3$, to
a quasihomogeneous one.

\par \vskip 1pc

Note that both the question and the conjecture above
are confirmed for $n = 2$.  Indeed, by the Lin--Zaidenberg
Theorem [LiZa] any contractible curve in  ${\bold C}^2$  is
equivalent to a quasihomogeneous one. But we do not know whether or not an
analogous statement holds in positive characteristic. So, we propose the
following

\par \vskip 1pc

{\bf PROBLEM.} {\it Given a field $k$ of positive characteristic, classify all
possible regular injections of the affine line ${\bold A}^1_k$ into the affine
plane ${\bold A}^2_k$.}

\par \vskip 1pc

Even for the embeddings this problem is wide open in char$k > 0$; see [D] and
references therein (this remark is due to P. Russell). Of course, in
characteristic $0$ this is the Epimorphism Theorem of Abhynkar and Moh and
Suzuki.

\par \vskip 1pc

{\bf 2.4.}  By  {\it an exotic} ${\bold C}^n$  we mean
a smooth affine algebraic variety diffeomorphic but
\par
non-isomorphic to  ${\bold C}^n$ (see [Za 2, 3]).

\par \vskip 1pc

{\bf QUESTION.}   Does there exist a pair of exotic ${\bold C}^n$'s
which are biholomorphic but not isomorphic ?
Does there exist a non-trivial deformation family of exotic
${\bold C}^n$'s  with the same underlying analytic structure ?

\par \vskip 1pc

{\bf 2.5.} {\bf CONJECTURE.}  Let  $X$  be a smooth algebraic variety
diffeomorphic to ${\bold C}^n$.  If $X$ is biholomorphic
to ${\bold C}^n$, then it is isomorphic to  ${\bold C}^n.$

\par \vskip 1pc

Note that by Ramanujam's Theorem [Ra] this is true for  $n = 2$.

\par \vskip 1pc

{\bf 2.6.}  {\bf PROBLEM.}  Verify that the hypersurface
$X = \{x + x^2y + z^3 + t^2 = 0 \} \subset {\bold C}^4$,
which is known to be an exotic  ${\bold C}^3$ (L. Makar-Limanov [ML; KML 2]),
is not biholomorphic to  ${\bold C}^3$.  The same for other known
hypersurfaces in  ${\bold C}^n$, which are exotic  ${\bold C}^{n-1}$'s (see [Za
3]).

\par \vskip 1pc

The next general problem is going back to F. Hirzebruch and
A. Van de Ven [VdV]; the only known positive result is the
Gurjar-Shastri Rationality Theorem for acyclic surfaces
(or, what is the same, for complex homology planes) [GuSh].

\par \vskip 1pc

{\bf 2.7.}  {\bf PROBLEM.}  Does there exist a non-rational exotic
${\bold C}^n$ ?

\par \vskip 1pc

{\bf 2.8.}  {\bf QUESTION.}   Let  $X$  be an exotic  ${\bold C}^n$.
Is it true that the action of the automorphism group
$\text{\rm Aut}\,X$  on  $X$  is not transitive ?
The same question for the group of analytic automorphisms.

\par \vskip 1.5pc

\S 3.  ELLIPTICALLY CONNECTED K\"AHLER MANIFOLDS

\par \vskip 1pc

{\bf 3.1.}  Note that an irreducible smooth quasiprojective
curve $C$  is non-hyperbolic iff the group  $\pi_1(C)$
is abelian.  Let  $X$  be a complex manifold.  We say that
$X$  is {\it elliptically connected} if any two points
$x,\,y \in X$  can be connected by a finite chain of holomorphic
images of non-hyperbolic curves and of complex tori.

\par \vskip 1pc

{\bf CONJECTURE.}  If  $X$  is an elliptically connected
complete K\"ahler manifold, then the group  $\pi_1(X)$
is almost abelian (or, in a weaker form, almost nilpotent).

\par \vskip 1pc

Here {\it `almost abelian'} resp. {\it `almost nilpotent'}
means that $\pi_1 (X)$ has an abelian resp.
nilpotent subgroup of a finite index.

\par \vskip 1pc

Note that {\it a rationally connected} (i.e. connected
by means of chains of rational curves) compact K\"ahler
manifold is simply connected (F. Campana [Ca]). Note also
that for quasiprojective surfaces the above conjecture
(in its stronger form) has been verified by using
the classification (L. Haddak--K. Oguiso--M. Zaidenberg). In higher dimensions,
it is unknown even whether or not $\pi_1 (X)$ is almost solvable.

On the other hand, it seems reasonable to enlarge the notion of an elliptically
connected manifold by admiting as members of chains the meromorphic images of
tori, or even of quasi--compact complex manifolds with abelian fundamental
groups.

And the last remark: the condition of k\"ahlerness is essential. As examples
one can consider quotients of $SL(2; {\bold C})$ by discrete cocompact
subgroups (J. Winkelmann; oral communication).

\par \vskip 1pc

REFERENCES :

[Ca] F. Campana, {\sl Remarques sur le rev\^etement universel des vari\'et\'es
K\"ahl\'eriennes compactes}, Bull. Soc. math. France {\bf 122} (1994), 255--284

[D] D. Daigle, {\sl A property of polynomial curves over a field of positive
characteristic}, Proc. Amer. Math. Soc.	{\bf 109} (1990), 887--894

[FlZa] H. Flenner, M. Zaidenberg, {\sl Q-acyclic surfaces and their
deformations}, Proc. Conf. "Classification of Algebraic Varieties", Mai 22--30,
1992, Univ. of l'Aquila, L'Aquila, Italy /Livorni ed. Contempor. Mathem. {\bf
162}, Providence, RI, 1994, 143--208

[Fu] T. Fujita, {\sl On the topology of non complete algebraic surfaces}, J.
Fac. Sci. Univ. Tokyo, Sect.IA, {\bf 29} (1982), 503--566

[GuMi] R.V. Gurjar, M. Miyanishi, {\sl Affine lines on logarithmic {\bf
Q}--homology planes}, Math. Ann. {\bf 294} (1992), 463--482

[GuPa] R.V. Gurjar, A.J. Parameswaran, {\sl Affine lines on {\bf Q}--homology
planes}, preprint, 1994, 1--18

[GuSh] R.V. Gurjar, A.R. Shastri,  {\sl On rationality of complex homology
2--cells}: I, II, J. Math. Soc. Japan {\bf 41} (1989), 37--56, 175--212

[KML 1] S. Kaliman, L. Makar-Limanov,  {\sl On morphisms into contractible
surfaces of Kodaira logarithmic dimension} 1,  preprint, 1994, 1--38

[KML 2] S. Kaliman, L. Makar-Limanov, {\sl On Russell's contractible
threefolds}, preprint,  1995, 1--22.

[LiZa] V. Lin, M. Zaidenberg, {\sl An irreducible simply connected curve
in ${\bold C}^{2}$ is equivalent to a quasihomogeneous curve}, Soviet Math.
Dokl., {\bf 28} (1983), 200-204

[ML] L. Makar-Limanov,  {\sl On the hypersurface $x+x^2y+z^2+t^3=0$
in $\bold C^4$}, preprint, 1994, 1--10

[MiTs] M. Miyanishi, S. Tsunoda, {\sl Absence of the affine lines on the
homology planes of general type}, J. Math. Kyoto Univ., {\bf 32} (1992),
443--450

[Ra] C.P. Ramanujam, {\sl A topological characterization of the affine
plane as an algebraic variety}, Ann. Math., 94 (1971), 69-88

[tD] T. tom Dieck,  {\sl Optimal rational curves and homology planes},
preprint, Mathematisches Institut, G\"ottingen, {\bf 9} (1992), 1--22

[VdV] A. Van de Ven,  {\sl Analytic compactifications of complex homology
cells}, Math. Ann. {\bf 147} (1962), 189--204

[Za 1] M. Zaidenberg, {\sl Isotrivial families of curves on affine surfaces and
characterization of the affine plane}, Math. USSR Izvestiya {\bf 30} (1988),
503-531. {\sl Addendum}, ibid {\bf 38} (1992), 435--437

[Za 2] M. Zaidenberg, An analytic cancellation theorem and exotic algebraic
structures on ${\bold C}^n$, $n \ge 3$,  Ast\'erisque {\bf 217} (1993),
251--282

[Za 3] M. Zaidenberg, {\sl On exotic algebraic structures on affine spaces},
preprint, Institute Fourier, Grenoble {\bf 305} (1995), 1--32

\enddocument